\documentclass[12pt]{article}
\usepackage[english]{babel}
\usepackage{amsfonts}
\usepackage[dvips]{graphicx}

\begin{document}

\title{\bf Optimal profiles in variable speed flows}

\date{February 2009}

\author{\bf Gianluca Argentini \\
\normalsize{[0,1]Bending - Italy}\\
\normalsize 01bending@gmail.com \\
\normalsize gianluca.argentini@gmail.com \\}

\maketitle

\begin{abstract}
Where a 2D problem of optimal profile in variable speed flow is resolved in a class of convex Bezier curves, using symbolic and numerical computations.\\

\end{abstract}

We study a particular problem on searching the optimal 2D shape or profile of a cartesian object immersed in a variable speed flow and subjected to some constraints on the boundary.\\

Let be $\{x,y\}$ a cartesian system, and $y = f(x)$ a function, defined for semplicity on the unit interval $[0,1]$, whose graph is the profile which we want to optimize. The function $f$ is subjected to the following constraints, arising from engineering requirements:\\

c1. $f(0) = 0$\\
\indent c2. $f(1) = 0$\\
\indent c3. $f$ is convex in $\left(0,1\right)$\\

\noindent Assume that the freestream flow is a vector field parallel to $y$-axis with a speed distribution according to a function $v=v(x)$ depending on $x$-variable. Also, assume that the flow pressure acting on the object profile is determined by Newton's sin-squared law (\cite{lachand}), that is, being $\rho$ the fluid constant density and $\alpha$ the angle between flow direction and tangent to profile,

\begin{equation}\label{sinSquared}
	p(x) = \rho v^2(x) sin^2\alpha = \rho v^2(x) \frac{1}{1+f'^2(x)}
\end{equation}

\noindent The local force acting on an element profile of $x$-extension $dx$ is $p(x)dx$, and this is the expression of the local resistance offered to the flow by the body. The total force on the profile is therefore

\begin{equation}\label{totalForce}
	F = \rho \int_0^1 \frac{v^2(x)}{1+f'^2(x)} dx
\end{equation}

\noindent The purpose of the mathematical analysis is the minimization of this functional in a suitable class of functions where one hopes to find the solution function $f(x)$. In general, usual methods of calculus of variations are applied, but one can easily test that technical difficulties arise in simple cases too. For example, let be $v(x)=ax$. The Euler-Lagrange equation (see \cite{lebedev}) associated to the minimization problem is than

\begin{equation}\label{EulerLagrange}
	0 = \frac{d}{dx}\frac{\partial F}{\partial f'} = - 2 \rho a^2 \frac{d}{dx}\frac{x^2 f'}{\left(1+f'^2\right)^2}
\end{equation}

\noindent from which a four-degree algebraic equations in the variable $f'$ arises, with $c$ arbitrary constant:

\begin{equation}
	x^2 f' + c \left(1+f'^2 \right)^2 = 0
\end{equation}

\noindent Such equation is not simple to manage and not so useful. It is more istructive, and useful for concrete applications, to consider a class of functions $f$ int the domain of Bezier curves.\\

Consider the three control points $\{\mathbf{P}_i, i=0,1,2\}$ = $\{0,0\}$, $\{a,1\}$,$\{1,0\}$, where $0 \leq a \leq 1$. The associated Bezier curve (see \cite{bartels}) is

\begin{equation}\label{BezierParam}
	b(t) = \left( b_1(t),b_2(t) \right) = \sum_{i=0}^2\mathbf{P}_i \left( _i^2 \right) t^i(1-t)^{2-i} = \{2at+(1-2a)t^2, 2t-2t^2\}
\end{equation}

\begin{figure}[h!]\label{BezierCurves}
	\begin{center}
	\includegraphics[width=5cm]{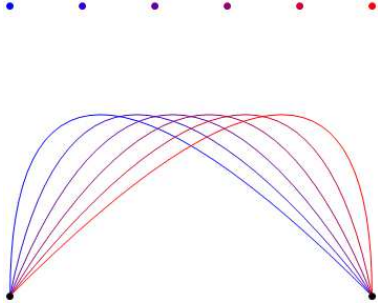}
	\caption{\small{\it Some Bezier curves in the case $a=0,0.2,0.4,0.6,0.8,1$, and relative control points.}}
	\end{center}
\end{figure}

\noindent Note that the parametric curve $b$ satisfies the conditions c1. and c2., respectively for $t=0$ and $t=1$, and c3. (see Fig.(1)).\\
\noindent Also, $x'(t) = 2a+2(1-2a)t$, and some usual algebraic computations show that $x'(t) \neq 0$ if $\frac{1}{2} < a < 1$. Under this assumption, if $f=f(x)$ is a cartesian representation of $b$, we have

\begin{equation}
	f'(x) = \frac{df}{dx} = \frac{d b_2}{dt} \frac{dt}{dx} = \frac{d b_2}{dt} \left[ \frac{d b_1}{dt} \right]^{-1} = \frac{1-2t}{a+(1-2a)t}
\end{equation}

\noindent and $dx = b_1'dt$, so that the functional (\ref{totalForce}) to minimize becomes

\begin{equation}\label{totalForceParam}
	F = F(a) = \rho \int_0^1 \frac{v^2(b_1)b_1'^3}{b_1'^2 + b_2'^2} dt
\end{equation}

\noindent Now the minimization problem is finding the value of the variable $a$, if exists, for which previous functional has its minimum value.\\

Consider, for example, the case of a freestream speed distribution $v(x) = -5x^3$; this situation has been proposed by a customer, as velocity field at oulet of a ventilating system and going towards a convex body. Along the profile of the Bezier curve (\ref{BezierParam}) the speed distribution has expression $v(x(t)) = -5 \left( 2at+(1-2a)t^2 \right)^3$ and the functional $F$ becomes

\begin{equation}\label{totalForceParamCustomer}
	F = F(a) = \rho \int_0^1 \frac{50t^6(a+t-2at)^3(2a+t-2at)^6}{(a+t-2at)^2+(1-2t)^2}dt
\end{equation}

\noindent We note that previous is an integral of a rational function and it is symbolically computable, but for our purpose we can only compute it numerically for a sufficient set of values of the parameter $a$, equally spaced in the interval $[0,1]$, and verify the existence of a local minimum for $F(a)$. Figure (2) shows the evidence of such a minimum for a value $0.6 < a_0 < 0.8$. An exact computation, by usual numerical techniques, gives the value $a_0 = 0.682564$ (our computation: use of {\bf FindMinimum}, by {\it Mathematica} v.7, (see \cite{wagon})).

\begin{figure}[h!]\label{Minimum}
	\begin{center}
	\includegraphics[width=6cm]{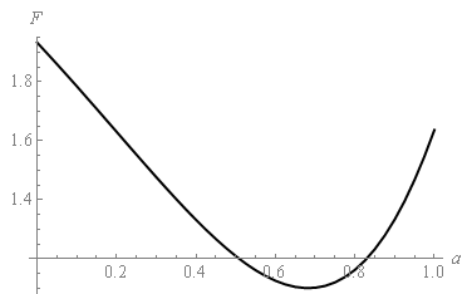}
	\caption{\small{\it Graph of $F(a)$ in the case $v(x)=-5x^3$.}}
	\end{center}
\end{figure}

\noindent The problem of finding the optimal profile in the suitable class of Bezier curves is so resolved in this concrete case.

\begin{figure}[h!]\label{profile}
	\begin{center}
	\includegraphics[width=5cm]{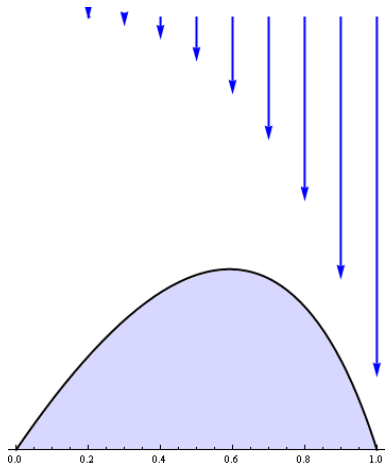}
	\caption{\small{\it The optimal profile. The point of maximum is in the region where the flow speed values are greater.}}
	\end{center}
\end{figure}

\begin{figure}[ht!]
	\includegraphics[width=2cm]{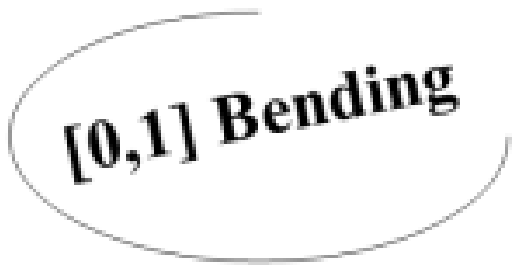}
\end{figure}
\noindent \tiny {\bf Gianluca Argentini}, mathematician, works on the field of\\
fluid dynamics, of acoustical noise reduction and optimization\\
of shapes for bodies moving inside fluid flows. He has found\\
{\bf [0,1]Bending}, a Design Studio in Italy dedicated to computational\\
engineering for scientific and industrial applications.

\end{document}